\documentclass[a4paper,12pt]{amsart}
\usepackage{amssymb}
\usepackage{ifthen}
\usepackage{graphicx}
\usepackage{float}

\newtheorem{thm}{Theorem}
\newtheorem{cor}{Corollary}
\newtheorem{lem}{Lemma}

\newtheorem{rem}{Remark}
\newtheorem{example}{Example}
\newtheorem{defn}{Definition}
\newtheorem{prob}{Problem}

\newtheorem{conj}{Conjecture}
\theoremstyle{definition}

\newcounter {own}
\def\theown {\thesection  .\arabic{own}}

\newenvironment{pf}[1][]{%
 \vskip 3mm
 \noindent
 \ifthenelse{\equal{#1}{}}%
  {{\slshape Proof. }}%
  {{\slshape #1.} }%
 }%
{\qed\bigskip}

\newcounter{alphabet}
\newcounter{tmp}
\newenvironment{Thm}[1][]{\refstepcounter{alphabet}%
\bigskip%
\noindent%
{\bf Theorem \Alph{alphabet}}%
\ifthenelse{\equal{#1}{}}{}{ (#1)}%
{\bf .} \itshape}{\vskip 8pt}

\newenvironment{Lem}[1][]{\refstepcounter{alphabet}%
\bigskip%
\noindent%
{\bf Lemma \Alph{alphabet}}%
{\bf .} \itshape}{\vskip 8pt}

\newcommand{\re}{{\rm Re\,}}

\newcommand{\IN}{{\mathbb N}}
\newcommand{\IC}{{\mathbb C}}
\newcommand{\ID}{{\mathbb D}}

\newcommand{\CC}{{\mathcal C}}

\def\be{\begin{equation}}
\def\ee{\end{equation}}

\newcommand{\bee}{\begin{enumerate}}
\newcommand{\eee}{\end{enumerate}}

\newcommand{\blem}{\begin{lem}}
\newcommand{\elem}{\end{lem}}
\newcommand{\bthm}{\begin{thm}}
\newcommand{\ethm}{\end{thm}}
\newcommand{\bcor}{\begin{cor}}
\newcommand{\ecor}{\end{cor}}
\newcommand{\beg}{\begin{example}}
\newcommand{\eeg}{\end{example}}
\newcommand{\begs}{\begin{examples}}
\newcommand{\eegs}{\end{examples}}
\newcommand{\bdefe}{\begin{defn}}
\newcommand{\edefe}{\end{defn}}
\newcommand{\bprob}{\begin{prob}}
\newcommand{\eprob}{\end{prob}}
\newcommand{\bei}{\begin{itemize}}
\newcommand{\eei}{\end{itemize}}

\newcommand{\bcon}{\begin{conj}}
\newcommand{\econ}{\end{conj}}
\newcommand{\bcons}{\begin{conjs}}
\newcommand{\econs}{\end{conjs}}
\newcommand{\bprop}{\begin{propo}}
\newcommand{\eprop}{\end{propo}}
\newcommand{\br}{\begin{rem}}
\newcommand{\er}{\end{rem}}
\newcommand{\brs}{\begin{rems}}
\newcommand{\ers}{\end{rems}}
\newcommand{\bo}{\begin{obser}}
\newcommand{\eo}{\end{obser}}
\newcommand{\bos}{\begin{obsers}}
\newcommand{\eos}{\end{obsers}}
\newcommand{\bpf}{\begin{pf}}
\newcommand{\epf}{\end{pf}}
\newcommand{\ba}{\begin{array}}
\newcommand{\ea}{\end{array}}
\newcommand{\beq}{\begin{eqnarray}}
\newcommand{\beqq}{\begin{eqnarray*}}
\newcommand{\eeq}{\end{eqnarray}}
\newcommand{\eeqq}{\end{eqnarray*}}

\newcommand{\ds}{\displaystyle}

\def\cc{\setcounter{equation}{0}   
\setcounter{figure}{0}\setcounter{table}{0}}

\def\cc{\setcounter{equation}{0}   
\setcounter{figure}{0}\setcounter{table}{0}}

\newcounter{minutes}
\divide\time by 60
\newcounter{hours}
\multiply\time by 60 \addtocounter{minutes}{-\time}

\begin{document}
\bibliographystyle{amsplain}
\title[
Harmonic Close-to-convex Functions and Minimal Surfaces ]{
Harmonic Close-to-convex Functions and Minimal Surfaces}

\thanks{
File:~\jobname .tex,
          printed: \number\year-\number\month-\number\day,
          \thehours.\ifnum\theminutes<10{0}\fi\theminutes}

\author{S. Ponnusamy 
}
\address{S. Ponnusamy, Department of Mathematics,
Indian Institute of Technology Madras, Chennai--600 036, India.}
\email{samy@iitm.ac.in}

\author{A. Rasila}
\address{A. Rasila, Department of Mathematics and Systems Analysis, Aalto University,
FI-00076 Aalto, Finland}
\email{antti.rasila@iki.fi}

\author{A. Sairam Kaliraj}
\address{A. Sairam Kaliraj, Department of Mathematics,
Indian Institute of Technology Madras, Chennai--600 036, India.}
\email{sairamkaliraj@gmail.com}



\begin{abstract}
In this paper, we study the family  $\CC_{H}^0$ of
sense-preserving complex-valued harmonic functions $f$ that are
normalized close-to-convex functions on the open unit disk $\ID$
with $f_{\overline{z}}(0)=0$. We derive a sufficient condition for
$f$ to belong to the class $\CC_{H}^0$. We take the analytic part of $f$ to
be $zF(a,b;c;z)$ or $zF(a,b;c;z^2)$  and for a suitable choice of
co-analytic part of $f$, the second complex dilatation
$w(z)=\overline{f_{\overline{z}}}/f_z$ turns out to be a square of
an analytic function. Hence $f$ is lifted to a minimal surface
expressed by an isothermal parameter. Explicit representation for
classes of minimal surfaces are given. Graphs generated by using
Mathematica are used for illustration.

\end{abstract}
\thanks{ }

\maketitle
\pagestyle{myheadings}
\markboth{S.Ponnusamy, A. Rasila, and A. Sairam Kaliraj}{On Harmonic Close-to-convex Functions and Minimal Surfaces}
\cc
\section{Introduction and Preliminary Results}

Denote by ${\mathcal H}$ the class of all complex-valued harmonic
functions $f$ in the unit disk ${\mathbb D}=\{z \in {\mathbb C}:\, |z|<1\}$ normalized by $f(0)=0=f_z(0)-1 $,
and let  ${\mathcal S}_H$ be the set of univalent functions in  ${\mathcal H}$. For
$f\in {\mathcal H}$, we have the canonical decomposition  $f=h+\overline{g}$, where $g$ and $h$ are analytic
on $\ID$. Here we call $h$ the analytic part of $f$ and $g$ the co-analytic part of $f$. We have
$$h(z)=z+\sum _{n=2}^{\infty}a_nz^n ~\mbox{ and }~g(z)=\sum _{n=1}^{\infty}b_nz^n, \quad z\in \ID
$$
and the Jacobian $J_f(z)$ of $f$ is
$$J_f(z)=|f_z(z)|^2-|f_{\overline z}(z)|^2= |h'(z)|^2-|g'(z)|^2.
$$
We say that the function $f$ is sense-preserving in $\ID$ if $J_f(z)>0$ in $\ID$.
According to a result of Lewy \cite{lewy-36}, the condition $J_f(z)>0$ in $\ID$ is necessary and sufficient for $f$
to be locally univalent and sense-preserving. Further,  $f\in {\mathcal H}$ is sense-preserving in $\ID$
if and only if $g'(z)=\omega (z)h'(z),$ where $\omega$ is analytic in $\ID$ with $|\omega (z)|<1$ in $\ID$.
We observe that if $g'(0)=0$, then $\omega$ fixes the origin so that by the Schwarz lemma one has $|\omega (z)|\leq |z|$
in $\ID$.

For basic results about the theory of planar harmonic mappings we refer to  \cite{Clunie-Small-84} and
the monograph of Duren \cite{Duren:Harmonic}. A function $f\in {\mathcal H}$ is said to be convex (starlike,
close-to-convex resp.) in $\ID_r:=\{z \in {\mathbb C}:\, |z|<r\}$ if it univalent in $\ID_r$ and $f(\ID_r)$
is convex (starlike with respect to the origin, close-to-convex resp.). By ${\mathcal K}_H$,  ${\mathcal S}_H^*$,
and ${\mathcal C}_H$, we denote the subclasses of functions in ${\mathcal S}_H$ that are convex, starlike, or
close-to-convex in the unit disk $\ID$, respectively. By ${\mathcal K}_H^0$,  ${\mathcal S}_H^{*0}$,
and ${\mathcal C}_H^0$, we mean the respective subclasses of functions $f$ from ${\mathcal K}_H$,  ${\mathcal S}_H^*$,
and ${\mathcal C}_H$ such that $f_{\overline {z}}(0)=b_1=0$.

One of the effective methods of constructing harmonic close-to-convex mappings from
conformal mappings is based on the following result due to  Clunie and Sheil-Small \cite{Clunie-Small-84}.

\begin{Lem}\label{uni-theo3}
If $h,g$ are analytic in $\mathbb D$ with $|h'(0)|>|g'(0)|$ and
$h+\epsilon g$ is close-to-convex for each $\epsilon$, $|\epsilon|=1$,
then $f=h+\overline{g}$ is close-to-convex in $\ID$.
\end{Lem}

This lemma  has been used to prove  many important results.
As a consequence of this result, in \cite{Hiroshi-Samy-2010},
the following result was established as a harmonic analog of Noshiro-Warschawski theorem
(see \cite[Theorem 2.16, p.~47]{Duren-book1}).

\begin{Lem}\label{uni-theo1}
{\rm \cite{Hiroshi-Samy-2010}}
Suppose $f=h+\overline{g}$ is harmonic on $\ID$ such that
${\rm Re\,} (e^{i\gamma}h'(z))>|g'(z)|$ for all $z\in \ID$, and for some $\gamma \in \mathbb{R}$.
Then $f$ is univalent, sense-preserving and close-to-convex in $\ID$.
\end{Lem}

In 1980, Mocanu \cite{Mocanu80} (see also \cite{Hiroshi-Samy-2010}) proved that if  $f=h+\overline{g}$ is
a harmonic mapping in a convex domain $\Omega$ such that ${\rm Re\,} (h'(z))>|g'(z)|$ for all $z\in \Omega$,
then $f$ is univalent and sense-preserving in $\Omega$.
An improved version of this results was given in \cite{KalajSamyMatti11,Hiroshi-Samy-2010}.
In order to discuss a general situation, it is appropriate to recall the following result due to Mocanu \cite{Mocanu80}.

\begin{Lem}\label{n2-theorem-A}
Let $G\in C^1(\ID)$ be univalent such that $G(\ID)$ is a convex domain and $J_G(z)>0$ in $\ID$.
Suppose that $F\in C^1(\ID)$, and
$$\re I(F,\overline{G})> |I(F,G)| ~\mbox{ for }z\in \ID,
$$
where
\[I(F,G)= \left|\begin{array}{cc} F_z & F_{\overline{z}} \\ G_z & G_{\overline{z}}\end{array}\right|.
\]
Then $F$ is sense-preserving and univalent in $\ID$.
\end{Lem}

Many functions that can be proved to be univalent by using
this lemma are also found to be close-to-convex in $\mathbb{D}$. In
\cite{PonSai-1(11)}, the following lemma was proved which, for example, leads to
the study of family of functions close-to-convex in $\mathbb{D}$
with respect to the convex function $-\log(1-z)$.

\begin{Lem} \label{APS1-11-lem1}
Let $f=h+\overline{g}$, where $f$ and $g$ are analytic in $\ID$ such that $h(0) = g(0) = 0$ and $h'(0) = 1$.
Further, let $G$ be univalent, analytic and convex in $\ID$. If $f$ satisfies
\be\label{APS1-11-eq1a}
{\rm Re}\left (e^{i\theta}\frac{h'(z)} {G'(z)}\right )>\left|\frac{g'(z)}{G'(z)}\right| ~\mbox{ for all $z\in \ID$ and
for some $\theta$ real},
\ee
then $f$ is sense-preserving, harmonic, univalent and close-to-convex in $\ID$.
\end{Lem}

The choice $G(z)=-\log(1-z)$ leads to the family
$${\mathcal F}=\{f\in {\mathcal H} :\, {\rm Re\,} \{(1-z)f_z (z)\}> |(1-z)f_{\overline{z}} (z)|, ~z\in\ID \}.
$$
According to Lemma \ref{APS1-11-lem1}, functions in ${\mathcal F}$
are close-to-convex in $\ID$. As a consequence of this result, the
following result  was established in \cite{PonSai-1(11)} together
with some applications associated with Gaussian hypergeometric functions.

\begin{Lem} \label{APS1-11-lem2} Suppose that $f=h+\overline{g} \in
{\mathcal H}$ satisfies the following condition
\be\label{APS1-11-eq2}
\sum _{n=1}^{\infty}|(n+1)a_{n+1}-n a_n | + \sum _{n=1}^{\infty}|(n+1)b_{n+1}-n b_n| \leq 1 - |b_1|
\ee
$(a_1=1)$. Then $f\in {\mathcal F}$. In particular, $f$ is
harmonic close-to-convex in $\ID$.
\end{Lem}

In fact one can obtain the following improved version of Lemma \ref{APS1-11-lem2}.

\bcor \label{APS1-11-cor1}
Suppose that  $f=h+\overline{g} \in {\mathcal H}$ satisfies the condition \eqref{APS1-11-eq2}.
Then $f\in {\mathcal F}_1$, where
$${\mathcal F}_1=\{f\in {\mathcal H} :\, |(1-z)h'(z) -1|<1- |(1-z)g'(z)|, ~z\in\ID \},
$$
and ${\mathcal F}_1\subset {\mathcal F}$.
\ecor

If we choose $G(z)=(1/2)\log((1+z)/(1-z))$ in Lemma \ref{APS1-11-lem1}, it leads to the family
$${\mathcal F}_2=\{f\in {\mathcal H} :\, {\rm Re\,}\{(1-z^2)f_z (z)\}> |(1-z^2)f_{\overline{z}} (z)|, ~z\in\ID \}.
$$
According to Lemma \ref{APS1-11-lem1}, functions in ${\mathcal F}_2$ are harmonic and close-to-convex in $\ID$.
To state a stronger version of the conclusion, it is necessary to recall the following definitions.

\bdefe\label{def2.1}
A domain $D\subset\mathbb{C}$ is called convex
in the direction $\alpha$ $(0\leq \alpha< \pi)$ if every line
parallel to the line through $0$ and $e^{i\alpha}$ has a connected
intersection with $D$.  A univalent harmonic function $f$ in
$\mathbb{D}$ is said to be {\it convex in the direction $\alpha$} if
$f(\mathbb{D})$ is convex in the direction $\alpha$.
\edefe

Obviously, every function that is convex in the direction $\alpha$
$(0\leq \alpha< \pi)$ is necessarily close-to-convex, but the converse is
not true.  The class of functions convex in one direction has been studied
by many mathematicians (see, for example, \cite{Do,Hen-Sch70,RZ})
as a subclass of functions introduced by Robertson \cite{Rober36}.
Now, we recall the following well-known result \cite{Clunie-Small-84}.

\begin{Lem}\label{lem2.1}
A harmonic function $f=h+\overline{g}$ locally univalent in $\mathbb{D}$ is a
univalent mapping of $\mathbb{D}$ onto a domain convex in the
direction of the real axis (resp. in the direction of the imaginary axis)
if and only if $h-g$ (resp. $h+g$) is a conformal
univalent mapping of $\mathbb{D}$ onto a domain convex  in the
direction of the real axis (resp. in the direction of the imaginary axis).
\end{Lem}

Paul Greiner \cite{Gre} has constructed numerous examples using the
method of shearing. However, by using this lemma, we can easily see that
functions in ${\mathcal F}_2$ are not only close-to-convex but is also convex
in the vertical direction.

\blem \label{PS2-12-lem1}
Functions in ${\mathcal F}_2$ are convex in the vertical direction.
\elem \bpf
Let $F=h+g$, where $f=h+\overline{g} \in {\mathcal F}_2$. Then
\beqq
{\rm Re\,}\{(1-z^2)F'(z)\}&=& {\rm Re\,}\{(1-z^2)(h'(z)+g'(z))\}\\
&\geq & {\rm Re\,}\{(1-z^2)h'(z)\}-|(1-z^2)g'(z)|\\
&=& {\rm Re\,}\{(1-z^2)f_z(z)\}-|(1-z^2)f_{\overline{z}}(z)| >0\cdot
\eeqq
From Theorem 1 of \cite{Hen-Sch70}, it is clear that $F$ is univalent and convex in the
vertical direction in $\mathbb{D}$. By \cite[Theorem 5.3]{Clunie-Small-84}, it is evident
that $f=h+\overline{g}$ is univalent and convex in vertical direction in $\mathbb{D}$.
\epf

 Now, we present a sufficient coefficient condition for functions to be in the family ${\mathcal F}_2$.

\blem \label{PS2-12-lem2}
Suppose that $f=h+\overline{g} \in
{\mathcal H}$ satisfies the following condition
\be\label{PS2-12-eq1}
\sum _{n=1}^{\infty}|(n+1)a_{n+1}-(n-1)a_{n-1} | + \sum _{n=1}^{\infty}|(n+1)b_{n+1}-(n-1) b_{n-1}| \leq 1 - |b_1|
\ee
$(a_1=1)$. Then $f\in {\mathcal F}_2$. In particular, $f$ is convex in the vertical direction in $\ID$ and
hence close-to-convex in $\ID$.
\elem
\begin{proof}
Without loss of generality, we may assume that $g(z) \not\equiv
0$. According to Lemma \ref{APS1-11-lem1}, it suffices to show
that \eqref{APS1-11-eq1a} holds for some convex function $G$. Now,
we set $G(z)=(1/2)\log((1+z)/(1-z))$. Then using \eqref{PS2-12-eq1} we find that
\beqq
{\rm Re\,}\left (\frac{f_z (z)} {G'(z)}\right )
&=& {\rm Re\,}\{(1-z^2)h'(z)\}\\
&=&{\rm Re\,}\left (1 + \sum _{n=1}^{\infty}\big ((n+1)a_{n+1}-(n-1) a_{n-1} \big )z^n  \right)\\
& \geq & 1-\sum _{n=1}^{\infty}|(n+1)a_{n+1}-(n-1) a_{n-1}|\\
& \geq &   |b_1| + \sum_{n=1}^{\infty}|(n+1)b_{n+1}-(n-1) b_{n-1}  |\\
& > & \left |b_1 + \sum_{n=1}^{\infty}\big ((n+1)b_{n+1}-(n-1)
b_{n-1}\big )z^n\right | \\
&=& |(1-z^2)g'(z)|= \left |\frac{f_{\overline{z}}(z)}{G'(z)}\right|.
\eeqq
The desired conclusion follows from Lemma \ref{APS1-11-lem1}.
\end{proof}

It is well known that the Euclidean coordinates of a minimal surface are harmonic functions of isothermal
parameters. The projection of such surface onto the base plane defines a harmonic mapping.
Conversely, a harmonic mapping which can be lifted to a minimal surface has a simple representation.
Weierstrass-Enneper representation (see \cite[p.177, Theorem]{Duren:Harmonic}) given below
describes the relation between a minimal surface defined by
isothermal parameters and the corresponding harmonic mappings.

\begin{Thm}\label{theorem_ms}
If a minimal graph $\left\{(u, v, F(u, v)): u+iv \in \Omega \right\}$
is parameterized by sense-preserving isothermal parameters $z=x+iy \in \mathbb{D}$, the projection
onto its base plane defines a harmonic mapping $w=u+iv=f(z)$ of $\mathbb{D}$ onto $\Omega$ whose
dilatation is the square of an analytic function. Conversely, if $f=h+\overline{g}$ is a
sense-preserving harmonic mapping of $\mathbb{D}$ onto some domain $\Omega$ with dilatation $w=q^2$
for some function $q$ analytic in $\mathbb{D}$, then the formulas
\be\label{ps2_ms}
u={\rm Re\,}(f(z)),~ v={\rm Im\,}(f(z)),~ t=2{\rm Im\,}\Big\{\int_0^z q(\zeta)h'(\zeta) \, \mathrm{d}\zeta \Big\}
\ee
define by isothermal parameters a minimal graph whose projection is $f$. Except for the choice of
sign and an arbitrary additive constant in the third coordinate function, this is the only such surface.
\end{Thm}


As an application,  we consider a particular type of analytic part of $f$, involving the Gaussian
hypergeometric function $F(a,b;c;z)$, and a suitable co-analytic part of $f$ so that dilatation
$\omega (z)$ turns out to be a constant multiple of $z^n$, $n \in \mathbb{N}$. When $n$ is an even
number, the dilatation is a square of an analytic function, hence the harmonic mapping can be
lifted to a minimal surface expressed by isothermal parameters.
Recently in \cite{PonQuRasila}, the authors considered certain class of
harmonic mappings convex in the horizontal direction with suitable dilatations $\omega$
and discussed the minimal surfaces associated with these harmonic mappings.


\section{Applications}

For complex numbers $a, b$ and $c$ with $ c\neq 0,-1,-2, \ldots $, the \textit{Gaussian} hyper-
geometric function defined by the series
$$F(a, b; c; z) = {}_2F_1(a,b;c;z) = \sum_{n=0}^{\infty}\frac{(a,n)(b,n)}{(c,n)(1,n)}z^n
$$
is analytic in $|z|<1$, where $(a, 0) = 1$ for $a \neq 0$, and
$(a, n) = a(a + 1)(a + 2) �\cdots (a +n-1)$ for $n\in \IN = \{1, 2, \ldots \}.$
For ${\rm Re}\, a > 0$ and  ${\rm Re}\,b > 0$,  we use the beta function $B(a,b)$ defined by
$$B(a, b) = \frac{\Gamma(a)\,\Gamma(b)}{\Gamma(a+b)},
$$
where $\Gamma (a)$ is the usual gamma function. In what follows, we need
the Stirling formula \cite[p.57, Equation (5)]{Bat53} given by
\be\label{APS1-11-eq3}
\lim_{n\rightarrow \infty} \frac{(a,n)(b,n)}{(c,n)(1,n)}=
\left \{ \begin{array}{lr}
\ds \frac {\Gamma(c)}{\Gamma(a)\Gamma(b)} & \mbox{if $c+1=a+b$},\\
0 & \mbox{if $c+1>a+b$},\\
\infty & \mbox{if $c+1<a+b$}.
\end{array}
\right.
\ee

\section{Main Results}

\bthm\label{PS2Ath1}
Let $a>0$, $b >0$, or  $a \in {\mathbb C}\setminus \{0\}$ with $b=\overline{a}$,
${\rm Re\,}a>0$.  Suppose that $a$, $b$, $m \in \mathbb{N}$ and $ \alpha \in {\mathbb C}$
are related by any one of the following conditions:
\be\label{ps2th1eq(a)}
ab \leq 1~\mbox{ and }~ |\alpha|(2B(a,b)-1)\leq 1,
\ee
\be\label{ps2th1eq(b)}
ab \geq \max\Big \{1, \frac{a+b}{2} \Big \}~\mbox{ and }~ |\alpha| \leq 2B(a,b)-1.
\ee
Then the harmonic function $f$ given by
\be\label{ps2th1eq(c)}
f(z)=zF(a,b;a+b;z)+ \ds \overline{\frac{\alpha z^{m+1}}{m+1}\Big\{ F(a,b;a+b;z)\ast F(2,m+1;m+2;z)\Big\} }
\ee
belongs to the class ${\mathcal F}_1$ $($and hence, is close-to-convex in $\ID$ with respect to
$-\log(1-z))$. Moreover, the dilatation of $f(z)$ is $\alpha z^m$.
\ethm
\begin{proof}
Following the standard notation, we let $f(z)=h(z)+\overline{g(z)}$, where
$$h(z)=zF(a,b;a+b;z)=\sum_{n=1}^{\infty}A_n z^n
$$
and
$$g(z)=\frac{\alpha z^{m+1}}{m+1}\Big\{ F(a,b;a+b;z)\ast F(2,m+1;m+2;z)\Big\}
=\frac{\alpha z^{m}}{m+1}\sum_{n=1}^{\infty}C_{n} z^{n},
$$
with
$$A_n=\frac{(a,n-1)(b,n-1)}{(a+b,n-1)(1,n-1)}, \quad n\geq 1,
$$
and
$$
C_n = A_n\frac{(2,n-1)(m+1,n-1)}{(m+2,n-1)(1,n-1)}, \quad n\geq 1 .
$$
Here $*$ denotes the usual Hadamard product (convolution) of power series.
We see that $A_1=1=C_1$, $A_n > 0$ and  $C_n >0 $ for all $n\geq 1$. Now, $f(z)$ takes the form
\beqq
f(z)=\sum_{n=1}^{\infty}A_n z^n + \frac{\overline{\alpha}}{m+1}\sum_{n=m+1}^{\infty}B_n \overline{z^n},
\eeqq
so that  $B_1=B_2=\cdots=B_m=0$ and $B_n=C_{n-m}$ for $n\geq m+1$. Further, a simple calculation gives
$$ nA_n-(n+1)A_{n+1}=\frac{A_n}{n(a+b+n-1)}X(n), \quad n\geq 1,
$$
where
$$X(n)=(n-1)(1-ab)+a+b-2ab.
$$
Similarly, a computation shows that for $n \geq m+1$,
$$nB_n-(n+1)B_{n+1}= \frac{nB_n}{(n-m)^2 (a+b+n-m-1)}Y(n),
$$
where
$$Y(n)=(n-m-1)(1-ab)+a+b-2ab.
$$
In order to apply Corollary \ref{APS1-11-cor1} together with \eqref{APS1-11-eq2}, it is convenient
to write
$$T:=T_1+(|\alpha|/(m+1))T_2,
$$
where
$$T_1=\sum _{n=1}^{\infty}|(n+1)A_{n+1}-n A_n | ~ \mbox{ and }~ T_2=\sum _{n=1}^{\infty}|(n+1)B_{n+1}-n B_n|,
$$
with $B_1=B_2=\cdots=B_m=0$. Clearly, by Corollary \ref{APS1-11-cor1}, $f\in {\mathcal F}_1$ 
(and hence, $f$ is close-to-convex with respect to $-\log(1-z)$) if $T \leq 1$.
Thus, to complete the proof, it suffices to show that $T \leq 1$ under the hypotheses of the theorem.

Case \textbf{(a)}: Suppose that $1-ab \geq 0$. Then
$$ \frac{1}{a}+\frac{1}{b} \geq \frac{1}{a}+a \geq 2,
$$
and so $a+b-2ab \geq 0$. In view of this observation and \eqref{ps2th1eq(a)}, it is clear that
$X(n)\ge X(1)\ge  0 $ for all $n\geq 1$. Similarly, $Y(n) \geq Y(m+1)\geq 0$ for $n \geq m+1$.
Thus, $T_1$ can be written as
\beqq
T_1&=& \lim_{k\rightarrow \infty}\sum _{n=1}^{k} \big (nA_n-(n+1)A_{n+1} \big )\\
&=& 1- \lim_{k\rightarrow \infty}(k+1)A_{k+1}\\
&=& 1- \lim_{k\rightarrow \infty}\left(k \frac{(a,k)(b,k)}{(a+b,k)(1,k)} + \frac{(a,k)(b,k)}{(a+b,k)(1,k)}\right)\\
&=& 1- \frac{ab}{a+b}  \lim_{k\rightarrow \infty} \frac{(a+1,k-1)(b+1,k-1)}{(a+b+1,k-1)(1,k-1)}
- \lim_{k\rightarrow \infty} \frac{(a,k)(b,k)}{(a+b,k)(1,k)}\cdot
\eeqq
Thus,  by the Stirling formula (\ref{APS1-11-eq3}), we get
\be\label{PS2-12-eq5}
T_1 = 1 - \frac{ab}{a+b} \frac{\Gamma(a+b+1)}{\Gamma(a+1) \Gamma(b+1)} = 1 - \frac {1}{B(a,b)}\cdot
\ee
Next, as $B_1=B_2=\cdots=B_m=0$ and $Y(n) \geq 0$ for all $n \geq 2$, we have
\beqq
T_2 &=& |(m+1)B_{m+1} - mB_m| + \sum _{n=m+1}^{\infty} \big |(n+1)B_{n+1} - nB_n \big |\\
&=& (m+1)+\lim_{k\rightarrow \infty}\sum _{n=m+1}^{k}\big (nB_n -(n+1)B_{n+1} \big )\\
&=& (m+1)+\lim_{k\rightarrow \infty}\big( (m+1)-(k+1)B_{k+1} \big )\\
&=& 2(m+1)-\lim_{k\rightarrow \infty}(k+1)B_{k+1} \cdot
\eeqq
In order to compute the limit on the right hand side, we rewrite $(k+1)B_{k+1}$ as
\beqq
(k+1)B_{k+1}&= &(k-m)\frac{(a,k-m)(b,k-m)}{(a+b,k-m)(1,k-m)}\frac{(2,k-m)(m+1,k-m)}{(m+2,k-m)(1,k-m)}\\
&& +(m+1)\frac{(a,k-m)(b,k-m)}{(a+b,k-m)(1,k-m)}\frac{(2,k-m)(m+1,k-m)}{(m+2,k-m)(1,k-m)}.
\eeqq
Applying the Stirling formula \eqref{APS1-11-eq3} as above, we easily obtain that
\be\label{PS2-12-eq6}
T_2 = 2(m+1)-\frac{m+1}{B(a,b)} \cdot
\ee
Combining \eqref{PS2-12-eq5} and \eqref{PS2-12-eq6}, we have
\beqq
T=1 + 2|\alpha|-\frac{1+|\alpha|}{B(a,b)}\cdot
\eeqq
Under the hypothesis \eqref{ps2th1eq(a)}, it is clear that $T \leq 1$ and therefore $f\in {\mathcal F}_1$.
\\

Case \textbf{(b)}: Suppose that \eqref{ps2th1eq(b)} holds. Then $ab \geq 1$ and
$2ab \geq a+b$ so that $ X(n)\leq  X(1)\leq 0$ for all $n\geq 1$. Similarly, it follows that
for all  $n \geq m+1$, $Y(n)\leq Y(m+1)\leq 0.$ Consequently,  as in the proof of Case \textbf{(a)},
the sum $T$ takes the form
\beqq
T&=& \lim_{k\rightarrow \infty}\sum _{n=1}^{k}\big ((n+1)A_{n+1}- nA_n \big) +
\frac{|\alpha|}{m+1} \lim_{k\rightarrow \infty}\sum _{n=m}^{k}\big((n+1)B_{n+1}- nB_n \big)\\
&=&\left (\frac{1}{B(a,b)}-1 \right) +  \frac{|\alpha|}{B(a,b)}.
\eeqq
Therefore $T\leq 1$ holds and thus, it follows that $f\in {\mathcal F}_1$
under the condition \eqref{ps2th1eq(b)} and the above relation.

Finally, from the power series representation of $h(z)$ and $g(z)$ it is easy to see that
$$g'(z)=\frac{\alpha}{m+1}\sum _{n=1}^{\infty}
\frac{(a,n-1)(b,n-1)}{(a+b,n-1)(1,n-1)}\frac{(2,n-1)(m+1,n-1)}{(m+2,n-1)(1,n-1)} (n+m)z^{n+m-1}.
$$
Since $(a,n)(a+n)=a(a+1,n)$, we may rewrite the last series as
$$g'(z)=\alpha z^m\sum _{n=1}^{\infty}\frac{(a,n-1)(b,n-1)}{(a+b,n-1)(1,n-1)} nz^{n-1}=\alpha z^m h'(z).
$$
Therefore the dilatation of $f(z)$ is $\omega (z)=g'(z)/h'(z)=\alpha z^m$.
\end{proof}
\bcor \label{PS2-12-cor1}
Assume the hypotheses of Theorem {\rm \ref{PS2Ath1}} on $a, b, m$ and $\alpha$. In addition, if
$m=2k$, then the formula $({\rm Re\,}(f(z)), {\rm Im\,}(f(z)), t(z))$
defines a minimal surface, where
$$f(z)=zF(a,b;a+b;z)+ \ds \overline{ \frac{\alpha z^{2k+1}}{2k+1}\Big\{ F(a,b;a+b;z)\ast F(2,2k+1;2k+2;z)\Big\} }
$$
and
$$t(z)=2{\rm Im\,}\Big\{\frac{\sqrt{\alpha} z^{k+1}}{k+1}[F(a,b;a+b;z)\ast F(2,k+1;k+2;z)]\Big\}+c,~ c \in \mathbb{R},
$$
whose projection is $f(z)$.
\ecor
\bpf
Theorem \ref{PS2Ath1} gives that $f(z)$ is univalent in $\mathbb{D}$. From the definition of $f(z)$, we see that
the second complex dilatation $\omega (z)=q^2(z)=\alpha z^{2k}$, which is a square of an analytic function.
Therefore by the Weierstrass-Enneper theorem (see Theorem \ref{theorem_ms}),
the function $f$ can be lifted to a minimal surface using the
formula \eqref{ps2_ms}. By considering the power series representation of $h(z)$ it is easy to see that $t(z)$
takes the form
$$t(z)=2{\rm Im\,}\Big\{\frac{\sqrt{\alpha} z^{k+1}}{k+1}[F(a,b;a+b;z)\ast F(2,k+1;k+2;z)]\Big\}+c,~ c \in \mathbb{R}.
$$
This completes the proof.
\epf

In the case $a,b>0$, we may reformulate Theorem \ref{PS2Ath1} in the following form.

\bcor \label{ps2-12-cor2}
Let $a>0, b>0, m\in \mathbb{N}$ and
$\alpha \in \mathbb{C}$ satisfies any one of the following
conditions
$$a \in (0, \infty),~b \in \Big(0, \frac{1}{a}\Big]~\mbox{ and }~ |\alpha|(2B(a,b)-1)\leq 1,
$$
$$a \in \Big(\frac{1}{2}, \infty \Big),~b \in \Big[\frac{a}{2a-1}, \infty
\Big)~\mbox{ and }~|\alpha|\leq 2B(a,b)-1.
$$
Then the harmonic function $f(z)$ defined in \eqref{ps2th1eq(c)} is close-to-convex in
$\mathbb{D}$.
\ecor


\beg
If we let $a=1$, and $b=1$ in Corollaries  \ref{PS2-12-cor1} and \ref{ps2-12-cor2},
then we have the following: If $m\in\mathbb{N}$ and $\alpha \in\IC$ such that $0<|\alpha|\leq 1$, then the function
$$f(z)=zF(1,1;2;z)+  \overline{\frac{\alpha z^{m}}{m+1} \int_0^z F(2,m+1;m+2;t)\,dt}
$$
belongs to ${\mathcal F}_1$. Using the derivative formula
$$F(a+1,b+1;c+1;z)=\frac{c}{ab}F'(a,b;c;z),
$$
the above integral can be computed and as a consequence,  we conclude that
$$f(z)=-\log(1-z)+   \overline{\frac{\alpha z^{m}}{m}(F(1,m;m+1;z)-1)}
$$
belongs to ${\mathcal F}_1$. In particular,
\be\label{ps2th1eq(f)}
f(z)= \left \{
\begin{array}{lr}
\ds -\log(1-z)-\overline{\alpha(z+\log(1-z))} & \mbox{if $m=1$},\\
-\log(1-z)-\overline{(\alpha/2)(2z+z^2+2\log(1-z))} & \mbox{if $m=2$},\\
-\log(1-z)-\overline{(\alpha/6)(6z+3z^2+2z^3+6\log(1-z))} & \mbox{if $m=3$},\\
-\log(1-z)-\overline{(\alpha/12)(12z+6z^2+4z^3+3z^4+12\log(1-z))} & \mbox{if $m=4$},
\end{array}
\right.
\ee
and
$$f(z)=-\log(1-z)-\overline{(\alpha/60)(60z+30z^2+20z^3+15z^4+12z^5+10z^6+60\log(1-z))}
$$
for $m=6$. Especially, when $m=2, 4, 6$ the above harmonic mappings can be lifted to minimal surface
in $\mathbb{R}^3$ whose coordinates are given by the formula $({\rm Re\,}(f(z)), {\rm Im\,}(f(z)), t(z))$ where
$$
t(z)=\left \{
\begin{array}{lr}
\ds -2{\rm Im\,}\{\sqrt{\alpha}(z + \log(1 - z))\} & \mbox{if $m=2$},\\
-{\rm Im\,}\{\sqrt{\alpha}(2 z + z^2 + 2 \log(1 - z))\} & \mbox{if $m=4$},\\
-(1/3){\rm Im\,}\{\sqrt{\alpha}(6 z + 3 z^2 +2 z^3 + 6 \log(1 - z))\} & \mbox{if $m=6$.}
\end{array}
\right.
$$
\eeg

\begin{figure}
\begin{center}
\includegraphics[height=6cm, width=5.5cm, scale=1]{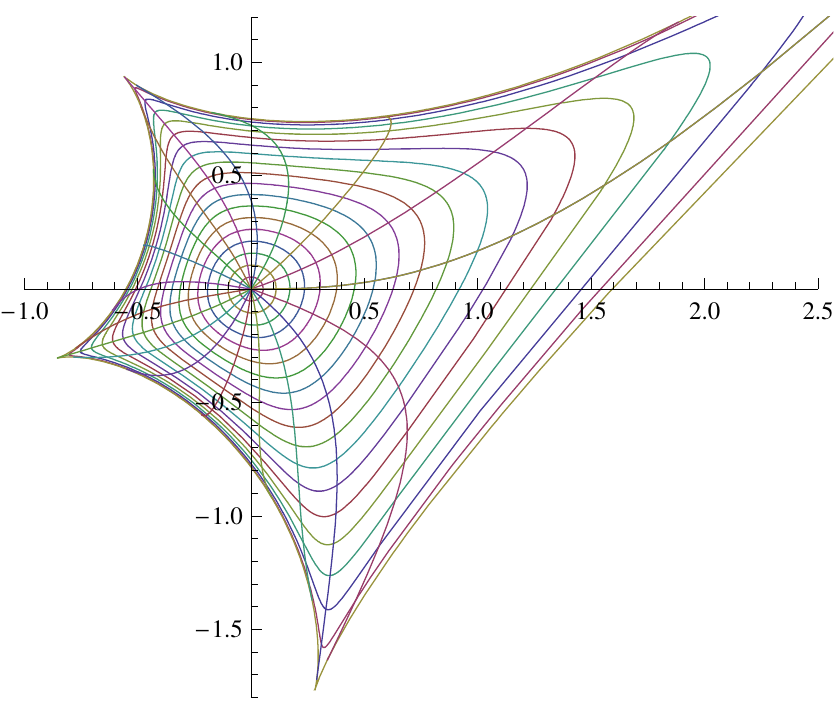}
\hspace{0.5cm}
\includegraphics[height=6cm, width=5.5cm, scale=1]{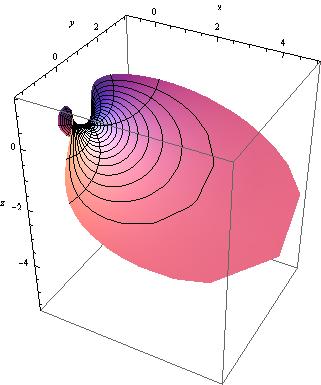}
\end{center}
\textbf{(a)} $m=2$ and $\alpha=-i$

\begin{center}
\includegraphics[height=6cm, width=5.5cm, scale=1]{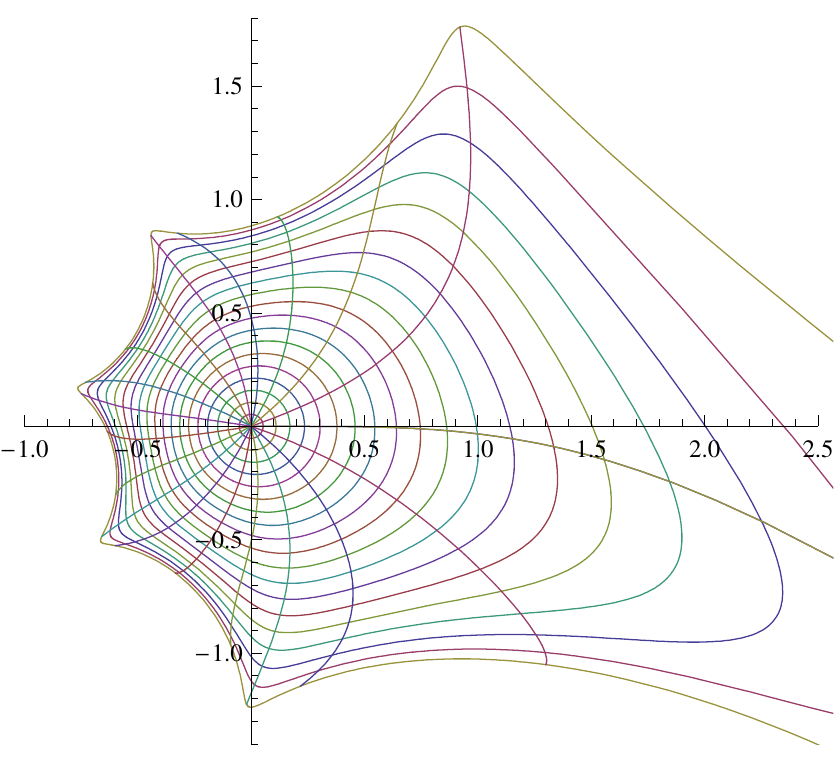}
\hspace{0.5cm}
\includegraphics[height=6cm, width=5.5cm, scale=1]{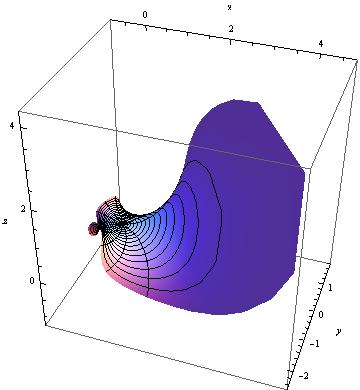}
\end{center}
\textbf{(b)} $m=4$ and $\alpha=0.75i$

\begin{center}
\includegraphics[height=6cm, width=5.5cm, scale=1]{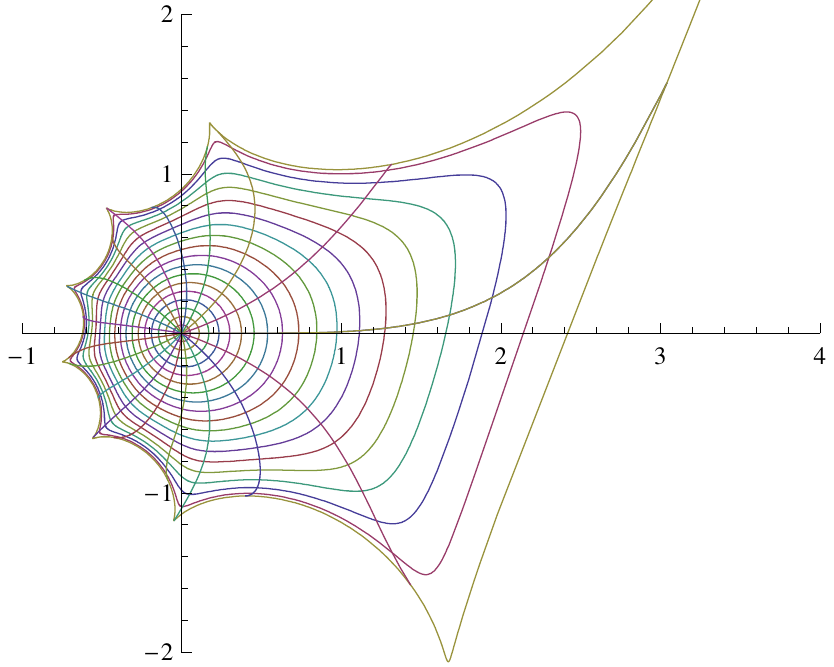}
\hspace{0.5cm}
\includegraphics[height=6cm, width=5.5cm, scale=1]{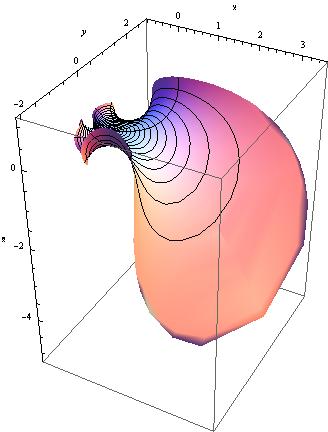}
\end{center}
\textbf{(c)} $m=6$ and $\alpha=-\frac{(1-i)}{\sqrt{2}}$
\caption{Images of $f(z)$ and the corresponding minimal surfaces for the indicated values of $m$ and $\alpha$ \label{pap2fig1}}

\end{figure}

The images of the disk $|z|<r$ for $r$ closer to $1$ under $f(z)$ in \eqref{ps2th1eq(f)}
for certain values of $m$, $\alpha$ and the corresponding minimal surfaces are shown in
Figures \ref{pap2fig1}\textbf{(a)}-\textbf{(c)}. These figures are drawn by using Mathematica.

%
\bthm\label{PS2Ath2}
Let $a>0$, $b >0$, or  $a \in {\mathbb C}\setminus \{0\}$ with $b=\overline{a}$,
${\rm Re\,}a>0$.  Suppose that $a$, $b$, $m \in \mathbb{N}$ and $ \alpha \in {\mathbb C}$
are related by any one of the following conditions:
\be\label{ps2th2eq(a)}
ab \leq \min\Big \{\frac{1}{2}, \frac{a+b}{3} \Big \}~\mbox{ and }~ |\alpha|(B(a,b)-1)\leq 1,
\ee
\be\label{ps2th2eq(b)}
ab \geq \max\Big \{\frac{1}{2}, \frac{a+b}{3}\Big \}~\mbox{ and }~ |\alpha| \leq B(a,b)-1.
\ee
Then the harmonic function $f$ given by
\be\label{ps2th2eq(c)}
f(z)=zF(a,b;a+b;z^2)+ \ds \overline{\frac{\alpha z^{m+1}}{m+1}\Big\{ F(a,b;a+b;z^2)\ast F(2,m+1;m+2;z)\Big\} }
\ee
belongs to ${\mathcal F}_2$ $($and hence, $f$ is convex in the vertical direction$)$.
Moreover, the dilatation of $f(z)$ is $\alpha z^m$.
\ethm
\begin{proof}
As in the proof of Theorem \ref{PS2Ath1}, we let $f(z)=h(z)+\overline{g(z)}$, where
$$h(z)=zF(a,b;a+b;z^2) = \sum_{n=1}^{\infty}A_{2n-1} z^{2n-1}
$$
and
$$g(z)=\frac{\alpha z^{m+1}}{m+1}\Big\{ F(a,b;a+b;z^2)\ast F(2,m+1;m+2;z)\Big\}
= \frac{\alpha z^m}{m+1}\sum_{n=1}^{\infty}C_{2n-1} z^{2n-1}
$$
with
$$A_{2n-1}=\frac{(a,n-1)(b,n-1)}{(a+b,n-1)(1,n-1)}, \quad n\geq 1,
$$
and
$$C_{2n-1}=A_{2n-1}\frac{(2,2n-2)(m+1,2n-2)}{(m+2,2n-2)(1,2n-2)}, \quad n\geq 1 .
$$
For convenience, we set $C_{2n-1}=B_{2n+m-1}$ so that
\beqq
f(z)=\sum_{n=1}^{\infty}A_{2n-1} z^{2n-1} + \overline{\frac{\alpha}{m+1}\sum_{n=1}^{\infty}B_{2n+m-1} z^{2n+m-1}}.
\eeqq
A simple calculation gives
$$ (2n-1)A_{2n-1}-(2n+1)A_{2n+1}=\frac{A_{2n-1}}{n(a+b+n-1)}X(n), \quad n\geq 1
$$
where
$$X(n)=(n-1)(1-2ab)+a+b-3ab.
$$
and by assumption \eqref{ps2th2eq(a)}, $X(n)\geq 0$ is satisfied. Similarly, for $n \geq 1$,
$$(2n+m-1)B_{2n+m-1}-(2n+m+1)B_{2n+m+1}= \frac{(2n+m-1)B_{2n+m-1}}{n(2n-1)(a+b+n-1)}X(n).
$$
Again, let $T:=T_1+(|\alpha|/(m+1))T_2$, where
$$T_1=\sum _{n=1}^{\infty}|(n+1)A_{n+1}-(n-1)A_{n-1} | ~ \mbox{ and }~ T_2=\sum _{n=1}^{\infty}|(n+1)B_{n+1}-(n-1)B_{n-1}|,
$$
with $B_1=B_2=\cdots=B_m=0$. By Lemma \ref{PS2-12-lem2}, it suffices to show that $T \leq 1$.

Following the proof of Theorem \ref{PS2Ath1}, the series $T_1$ and $T_2$ may be computed easily
and see that
\beqq
T=1 -\frac{2}{B(a,b)}+2|\alpha|-\frac{2|\alpha|}{B(a,b)}.
\eeqq
Under the hypothesis \eqref{ps2th2eq(a)},  $T \leq 1$ and therefore $f\in {\mathcal F}_2$.

Similarly, when \eqref{ps2th2eq(b)} holds, $X(n) \leq 0$ so that $T$ reduces to
\beqq
T=\frac{2}{B(a,b)}-1+\frac{2|\alpha|}{B(a,b)}.
\eeqq
Therefore the conclusion $T\leq 1$ follows from \eqref{ps2th2eq(b)} and the above relation.

Finally, from the power series representation of $h(z)$ and $g(z)$ it is easy to see that
$$g'(z)=\frac{\alpha}{m+1}\sum _{n=1}^{\infty}
\frac{(a,n-1)(b,n-1)}{(a+b,n-1)(1,n-1)}\frac{(2,2n-2)(m+1,2n-2)}{(m+2,2n-2)(1,2n-2)} (2n+m-1)z^{2n+m-2}.
$$
Since $(a,n)(a+n)=a(a+1,n)$, we may rewrite the last series as
$$g'(z)=\alpha z^m\sum _{n=1}^{\infty}\frac{(a,n-1)(b,n-1)}{(a+b,n-1)(1,n-1)} (2n-1) z^{2n-2}=\alpha z^m h'(z).
$$
Therefore the dilatation of $f(z)$ is $\omega (z)=g'(z)/h'(z)=\alpha z^m$.
\end{proof}
\bcor \label{PS2-12-cor2}
Suppose $a, b, m$ and $\alpha$ satisfies the hypothesis of Theorem \ref{PS2Ath1} and further $m=2k$.
Then the formula $({\rm Re\,}(f(z)), {\rm Im\,}(f(z)), t(z))$ defines a minimal surface, where
$$f(z)=zF(a,b;a+b;z^2)+ \ds \overline{\frac{\alpha z^{2k+1}}{2k+1}\Big\{ F(a,b;a+b;z^2)\ast F(2,2k+1;2k+2;z)\Big\} }
$$
and
$$t(z)=2{\rm Im\,}\Big\{\frac{\sqrt{\alpha} z^{k+1}}{k+1}[ F(a,b;a+b;z^2)\ast F(2,k+1;k+2;z)]\Big\}+c,~ c \in \mathbb{R},
$$
whose projection is $f(z)$.
\ecor
\bpf The result follows from the
previous theorem and the Weierstrass-Enneper representation for
minimal surface whose coordinates $(u,v,t)$ in $\mathbb{R}^3$ is
given by
$$u={\rm Re\,}(f(z)),~ v={\rm Im\,}(f(z)),~ t(z)=2{\rm Im\,}\Big\{\int_0^z q(\zeta)h'(\zeta) \, \mathrm{d}\zeta \Big\}.
$$
By considering the power series representation of $h(z)$ it is easy to see that $t(z)$ takes the form
$$t(z)=2{\rm Im\,}\Big\{\frac{\sqrt{\alpha} z^{k+1}}{k+1}[ F(a,b;a+b;z^2)\ast F(2,k+1;k+2;z)]\Big\}+c,~ c \in \mathbb{R}.
$$
This completes the proof.
\epf

\br\label{PS2rm2}
{\rm From Lemma \ref{PS2-12-lem1} it is clear that the function $f(z)$ in \eqref{ps2th2eq(c)} is not
only close-to-convex in $\mathbb{D}$, but also convex in the vertical direction in $\mathbb{D}$. Using
\cite[Theorem 5.3]{Clunie-Small-84}, it is easy to see that the conformal pre-shear of $f(z)$ defined by
$$ \phi(z)= h(z)+g(z)=zF(a,b;a+b;z^2)+ \ds \frac{\alpha z^{m+1}}{m+1}\Big\{ F(a,b;a+b;z^2)\ast F(2,m+1;m+2;z)\Big\}
$$
is univalent and convex in the vertical direction in $\mathbb{D}$.
}
\er

\bcor\label{PS2cor3}
Let $a>0$, $b >0$ be such that
$$b\in \left \{ \begin{array}{rl}
\ds \left (0, \frac{1}{2a}\right ] & \mbox{ if $a\in (0,\frac{1}{2}]\cup [1,\infty)$},\\[4mm]
\ds \left (0, \frac{a}{3a-1}\right ] & \mbox{ if $a\in [\frac{1}{2},1]$,}
\end{array}
\right .
$$
and $\alpha \in \IC$ satisfies the condition $|\alpha|(B(a,b)-1) \leq 1$. Then the function $f(z)$ defined in
\eqref{ps2th2eq(c)} belongs to the class ${\mathcal F}_2$.
\ecor\bpf
In order to prove the result, it suffices to show that the above
conditions imply \eqref{ps2th2eq(a)}. When $0< b \leq (1/2a)$, it
is clear that $ab \leq 1/2$. If $0<a<1/3$, then $3a-1<0$ and hence
$3ab<a+b$. We have
$$ b \leq \frac{1}{2a} \leq \frac{a}{3a-1}, ~\mbox{ if }~ a \in \Big(\frac{1}{3}, \frac{1}{2}\Big]\cup [1, \infty),
$$
which again implies $3ab \leq a+b$. When a=1/3, the inequality
$3ab \leq a+b$ holds.
$$ b \leq \frac{a}{3a-1} \leq \frac{1}{2a} ~\mbox{ if }~ a \in \Big[\frac{1}{2}, 1\Big].
$$
Combining these observations, we get the result.
\epf


\bcor\label{PS2cor4}
Let $a>0$, $b >0$ be such that
$$b\in \left \{ \begin{array}{rl}
\ds \left [\frac {a}{3a-1}, \infty\right ) & \mbox{ if $a\in (\frac{1}{3},\frac{1}{2}]\cup [1,\infty)$},\\[4mm]
\ds \left [\frac{1}{2a}, \infty\right ) & \mbox{ if $a\in [\frac{1}{2},1]$},
\end{array}
\right.
$$
and $ \alpha \in \IC$ satisfies the condition
$|\alpha| \leq B(a,b)-1$. Then the function $f(z)$ defined in
\eqref{ps2th2eq(c)} belongs to the class ${\mathcal F}_2$.
\ecor \bpf
The condition
$$ ab \geq \max\Big \{\frac{1}{2}, \frac{a+b}{3}\Big \}
$$
can be rewritten as
$$ b\geq \max \left \{ \frac {1}{2a}, \frac {a}{3a-1}\right \} =
\left \{ \begin{array}{lr}
\ds \frac {a}{3a-1} & \mbox{ for $(a-1)(a-1/2)\geq 0$},\\[4mm]
\ds \frac {1}{2a} & \mbox{ for $(a-1)(a-1/2)\leq 0$}.
\end{array}
\right.
$$
Since $a$ is positive, $a/(3a-1)>0$ if $a>1/3$. The result follows from the above inequality and
\eqref{ps2th2eq(b)}.
\epf

\beg
If $b \in(0, 1/2)$, $m\in \mathbb{N}$ and $\alpha \in \IC$ such that $0<|\alpha|\leq b/(1-b)$, then the function
$$f(z)=zF(1,b;1+b;z^2)+ \ds \overline{\frac{\alpha z^{m+1}}{m+1}\Big\{ F(1,b;1+b;z^2)\ast F(2,m+1;m+2;z)\Big\} }
$$
belongs to the class ${\mathcal F}_2$ and hence $f$ is close-to-convex in $\ID$.
\eeg

\begin{figure}
\begin{center}
\includegraphics[height=6cm, width=5.5cm, scale=1]{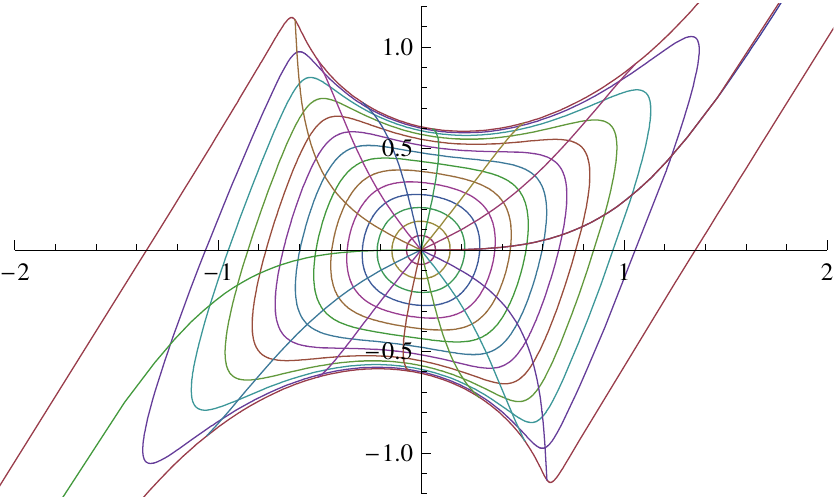}
\hspace{0.5cm}
\includegraphics[height=6cm, width=5.5cm, scale=1]{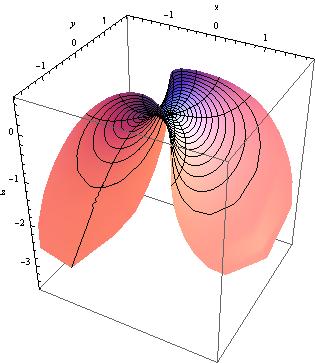}
\end{center}
\textbf{(e)} $a=\frac{3}{4}, b=\frac{2}{3}, m=2$ and $\alpha=-\frac{0.871(1+i)}{\sqrt{2}}$
\caption{Images of $f(z)$ and the corresponding minimal surfaces for the indicated values of $a, b, m$ and $\alpha$ \label{pap2fig2}}
\end{figure}

The images of the disk $|z|<r$ for $r$ closer to 1 under $f(z)$ in \eqref{ps2th2eq(c)}
for certain values of $a, b, m$, $\alpha$ and the corresponding minimal surfaces are shown in
Figures \ref{pap2fig2}\textbf{(a)}-\textbf{(e)}. As mentioned in the Remark \ref{PS2rm2} functions present in the class
${\mathcal F}_2$ are convex in the vertical direction.

\begin{figure}[H]
\begin{center}
\includegraphics[height=5.5cm, width=5cm, scale=1]{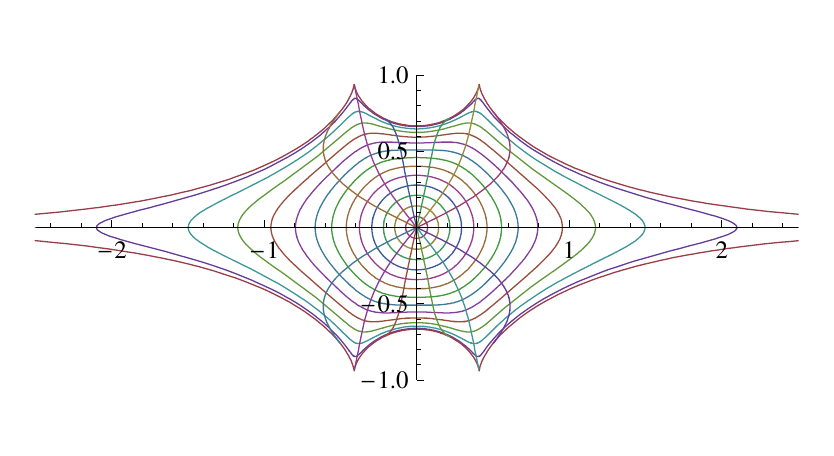}
\hspace{0.5cm}
\includegraphics[height=5.5cm, width=5cm, scale=1]{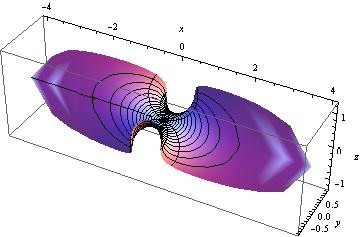}
\end{center}
\textbf{(a)} $a=1, b=0.5, m=4$ and $\alpha=1$
\end{figure}
\begin{figure}[H]
\begin{center}
\includegraphics[height=5.5cm, width=5cm, scale=1]{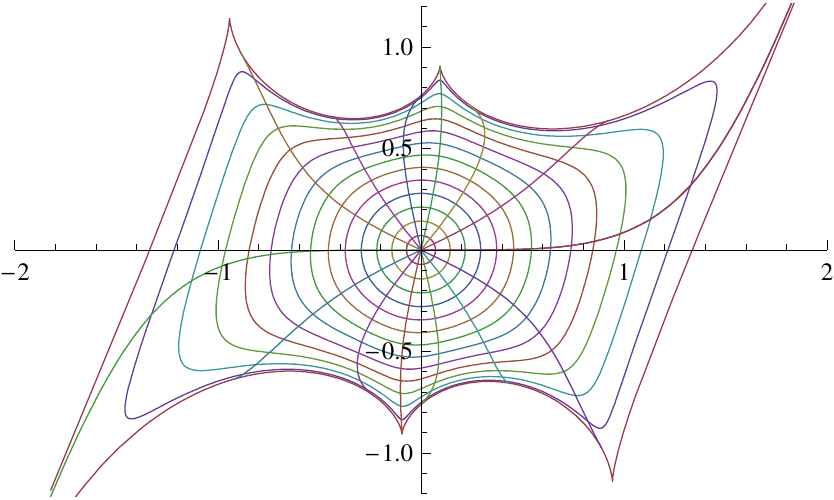}
\hspace{0.5cm}
\includegraphics[height=5.5cm, width=5cm, scale=1]{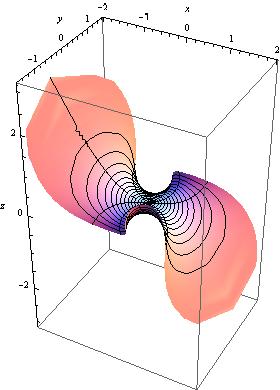}
\end{center}
\textbf{(b)} $a=1, b=0.5, m=4$ and $\alpha=-\frac{(1+i)}{\sqrt{2}}$
\begin{center}
\includegraphics[height=6cm, width=5.5cm, scale=1]{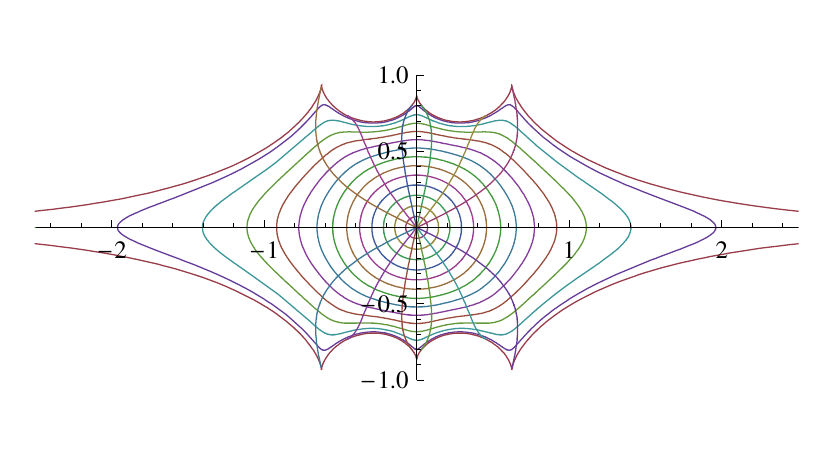}
\hspace{0.5cm}
\includegraphics[height=6cm, width=5.5cm, scale=1]{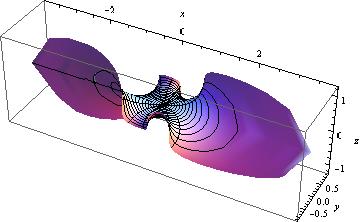}
\end{center}
\textbf{(c)} $a=1, b=0.5, m=6$ and $\alpha=1$
\begin{center}
\includegraphics[height=6cm, width=5.5cm, scale=1]{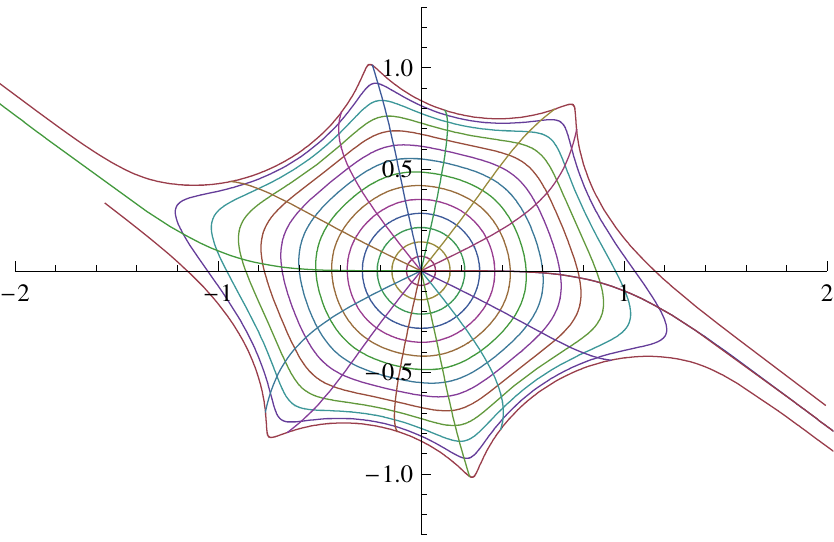}
\hspace{0.5cm}
\includegraphics[height=6cm, width=5.5cm, scale=1]{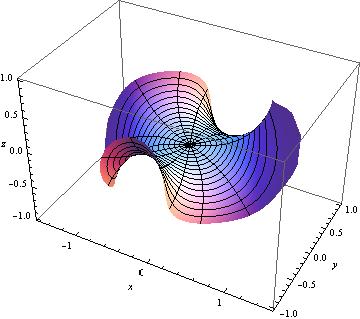}
\end{center}
\textbf{(d)} $a=1, b=\frac{1}{7}, m=4$ and $\alpha=\frac{i}{7}$
\end{figure}
\begin{figure}[H]
\begin{center}
\includegraphics[height=6cm, width=5.5cm, scale=1]{T2HM_1_0_5_6_1.pdf}
\hspace{0.5cm}
\includegraphics[height=6cm, width=5.5cm, scale=1]{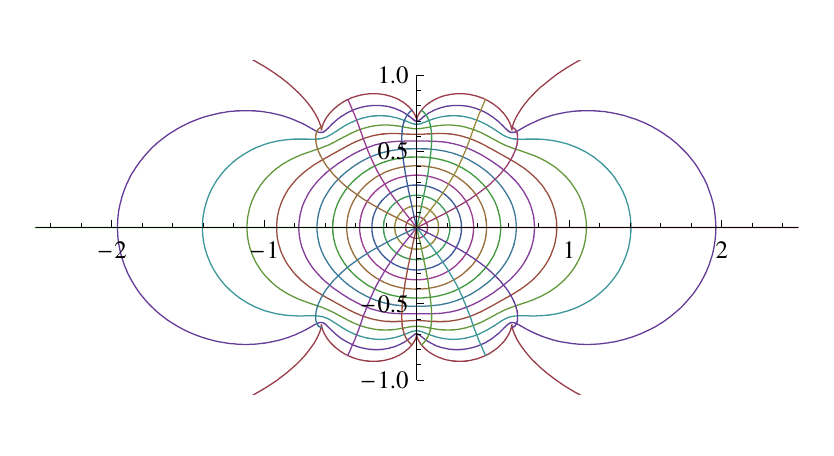}
\end{center}
\textbf{(c)} $a=1, b=0.5, m=6$ and $\alpha=1$
%
\begin{center}
\includegraphics[height=6cm, width=5.5cm, scale=1]{T2HM_1_1O7_4_3O4_piby2.pdf}
\hspace{0.5cm}
\includegraphics[height=6cm, width=5.5cm, scale=1]{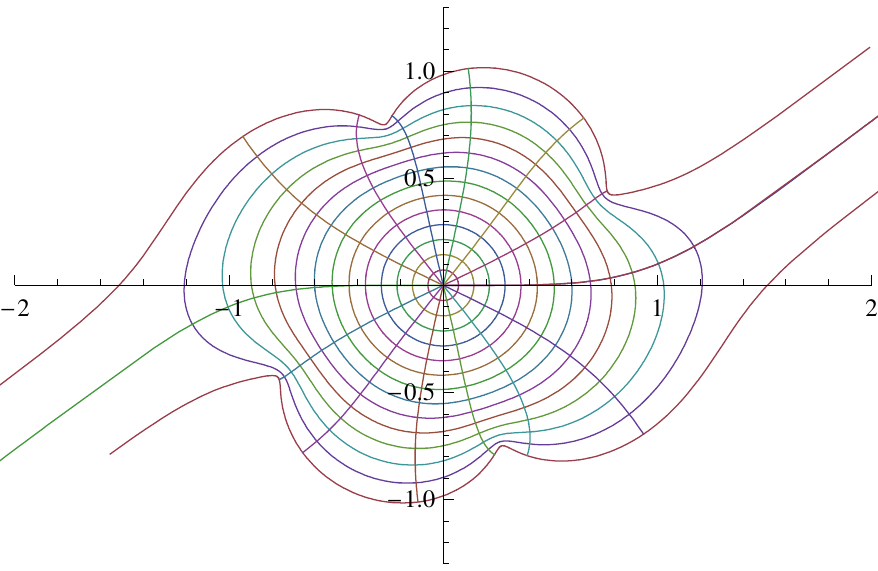}
\end{center}
\textbf{(d)} $a=1, b=\frac{1}{7}, m=4$ and $\alpha=\frac{i}{7}$

\begin{center}
\includegraphics[height=6cm, width=5.5cm, scale=1]{T2HM_3O4_2O3_2_p871_5piby4.pdf}
\hspace{0.5cm}
\includegraphics[height=6cm, width=5.5cm, scale=1]{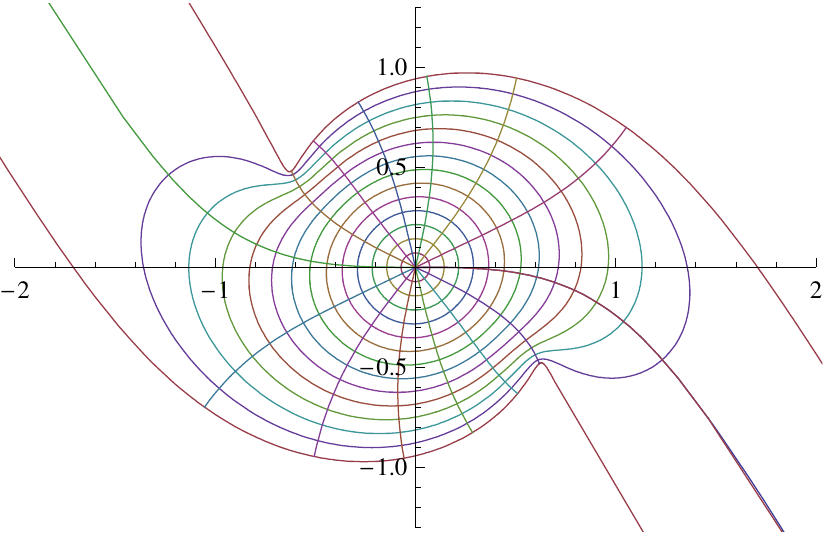}
\end{center}
\textbf{(e)} $a=\frac{3}{4}, b=\frac{2}{3}, m=2$ and $\alpha=-\frac{0.871(1+i)}{\sqrt{2}}$
\caption{Images of $f(z)$ and the corresponding conformal pre-shears for the indicated
values of $a, b, m$ and $\alpha$ \label{pap2fig3}}

\end{figure}

The image of the harmonic mappings in
Figures \ref{pap2fig2}\textbf{(a)}-\textbf{(e)} together with the corresponding conformal pre-shears are shown
in Figures \ref{pap2fig3}\textbf{(a)}-\textbf{(e)}.
\begin{figure}[H]
\begin{center}
\includegraphics[height=6cm, width=5.5cm, scale=1]{T2HM_1_0_5_4_1.pdf}
\hspace{0.5cm}
\includegraphics[height=6cm, width=5.5cm, scale=1]{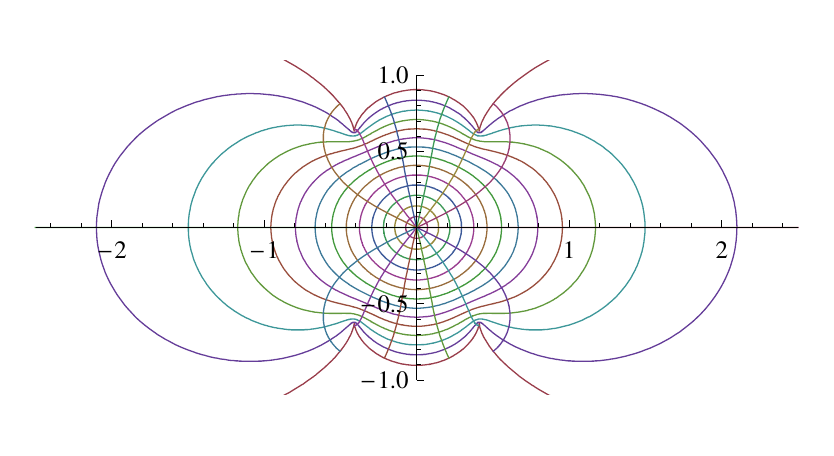}
\end{center}
\textbf{(a)} $a=1, b=0.5, m=4$ and $\alpha=1$
\begin{center}
\includegraphics[height=6cm, width=5.5cm, scale=1]{T2HM_1_0_5_4_5piby4.pdf}
\hspace{0.5cm}
\includegraphics[height=6cm, width=5.5cm, scale=1]{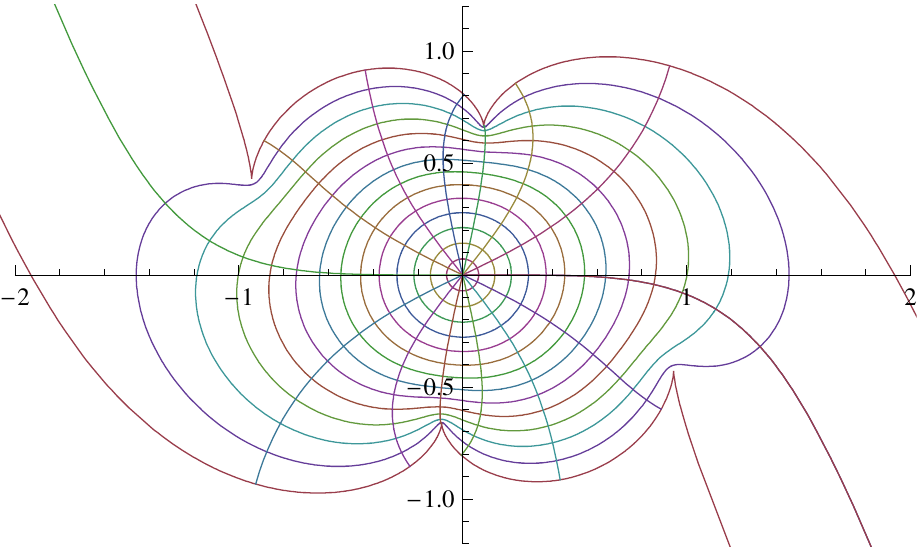}
\end{center}
\textbf{(b)} $a=1, b=0.5, m=4$ and $\alpha=-\frac{(1+i)}{\sqrt{2}}$
\end{figure}


\begin{thebibliography}{150}

\bibitem{Bat53} H.~Bateman~(Ed. by A.~Erdelyi, W.~Magnus, F.~Oberhettinger, and F.~G.~Tricomi),
Higher Transcendental Functions, Vol. I, McGraw-Hill, New York, 1953.



%

\bibitem{Clunie-Small-84} J.~G.~Clunie and T.~Sheil-Small,
\emph{Harmonic univalent functions},
Ann. Acad. Sci. Fenn. Ser. A.I. {\bf 9}(1984), 3--25.

\bibitem{Do} M.~Dorff,
Convolutions of planar harmonic convex mappings,
\textit{Comp. Vari. Theo. Appl.} {\bf 45}(2001), 263--271.

\bibitem {Duren-book1} P.~Duren, Univalent Functions
(Grundlehren der mathematischen Wissenschaften 259, New York, Berlin, Heidelberg, Tokyo),
Springer-Verlag, 1983.

\bibitem{Duren:Harmonic} P.~Duren,
\textrm{Harmonic Mappings in the Plane},
Cambridge Tracts in Mathematics,
\textbf{156}, Cambridge Univ. Press, Cambridge, 2004.

\bibitem{Gre} P. Greiner,
Geometric properties of harmonic shears,
\textit{Comput. Methods Funct. Theory} {\bf 4}(1)(2004), 77--96.

\bibitem{Hen-Sch70} W.~Hengartner, and G.~Schober,
\emph{On Schlicht Mappings to Domains Convex in One Direction},
\textrm{Comment. Math. Helv.} {\bf 45} (1970), 303-314.



\bibitem{KalajSamyMatti11} D.~Kalaj, S.~Ponnusamy, and M.~ Vuorinen,
\emph{Radius of close-to-convexity of harmonic functions},
\textit{Complex Var. Elliptic Equ.}, Revised.

\bibitem{lewy-36} H.~Lewy,
\emph{On the nonvanishing of the Jacobian in certain one-to-one mappings},
\textrm{Bull. Amer. Math. Soc.} {\bf 42} (1936), 689--692.

\bibitem{Mocanu80} P.~T.~Mocanu,
\emph{Sufficient conditions of univalency for complex functions in the class $C^{1}$},
Anal. Num\'er. Th\'eor. Approx. \textbf{10}(1)(1981),  75--79.


%


\bibitem{PonQuRasila}
 S.~Ponnusamy, T.~Quach, and A.~Rasila,
\emph{Harmonic shears of slit and polygonal mappings},
 {\tt  http://arxiv.org/pdf/1201.2015.pdf}

\bibitem {PonSai-1(11)}  S.~Ponnusamy, and A.~Sairam Kaliraj
\emph{On Harmonic Close-to-convex Functions},
Comput. Methods Funct. Theory, (2012), To appear.

\bibitem {Hiroshi-Samy-2010}  S.~Ponnusamy, H.~Yamamoto and H.~Yanagihara,
\emph{Variability regions for certain families of harmonic univalent mappings},
Complex Var. Elliptic Equ. (2011), in print.

\bibitem{Rober36} M. S. Robertson,
Analytic functions starlike in one direction, \textit{Am. J. Math.}
\textbf{58}(1936), 465--472.

\bibitem{RZ} W. C. Royster and M. Ziegler,
Univalent functions convex in one direction,
\textit{Publ. Math. Debrecen} \textbf{23}(1976), 339--345.







\end{thebibliography}
\end{document}